\documentclass[final]{siamltex}

\usepackage{amssymb,amsmath}
\usepackage{epsfig}
\usepackage{psfrag}

\newcommand{\tp}{^\mathrm{T}} 

\title{A dual eigenvector condition for strong lumpability of Markov chains\thanks{This work was funded by PACE (Programmable Artificial Cell Evolution), a European Integrated Project in the EU FP6-IST-FET Complex Systems Initiative, and by EMBIO (Emergent Organisation in Complex Biomolecular Systems), a European Project in the EU FP6-NEST Initiative.}}

\author{Martin Nilsson Jacobi\thanks{Complex Systems Group, Department of Energy and Environment, Chalmers University of Technology, 412 96 G\"{o}teborg, Sweden ({\tt mjacobi@chalmers.se} and {\tt olofgo@chalmers.se}).} 
        \and Olof G\"{o}rnerup\footnotemark[2]}
   
\begin{document}

\maketitle

\begin{abstract}
Necessary and sufficient conditions for identifying strong lumpability in Markov chains are presented. We show that the states in a lump necessarily correspond to identical elements in eigenvectors of the dual transition matrix. If there exist as many dual eigenvectors that respect the necessary condition as there are lumps in the aggregation, then the condition is also sufficient. The result is demonstrated with  two simple examples.
\end{abstract}

\begin{keywords} 
Strong lumpability; Aggregated Markov chains; State space reduction
\end{keywords}

\begin{AMS}
60J22, 15A18, 60J99
\end{AMS}

\pagestyle{myheadings}
\thispagestyle{plain}
\markboth{M. NILSSON JACOBI AND O. G\"{ORNERUP}}{STRONG LUMPABILITY OF MARKOV CHAINS}

\section{Introduction}
In this paper we present a method for coarse graining Markov chains by partitioning, or lumping, the state space into aggregated Markov chains. This approach to coarse grain Markov chains has been studied for a long time, see e.g.~\cite{Kemeny99,Rogers}. Due to the combinatorial explosion of the number of possible partitions of the state space, the conditions used to define lumpability are not useful for identifying aggregated states even for relatively moderately sized Markov chains. Existing efficient methods for lumping are concerned with the special  task of finding the most coarse approximate lumping  \cite{Derisavi03}, or special situations where the transition matrix possesses known symmetries (e.g. permutation symmetries)~\cite{Buchholz, Obal,Ring}. There also exist efficient methods for testing that a given partition constitute a valid lumping \cite{Barr77, Jernigan03,Gurvits}. It should be noted that there can exist multiple different accepted lumpings of a Markov chain. Furthermore, different lumpings are in general not related through successive refinements. This paper addresses the general question: Given the transition matrix of a Markov process, how can the state space be partitioned such that the Markov property is not lost in the aggregated process?

\section{Main result}
Consider a Markov chain $X_t$, $t=0,1,\dots$, with a finite state space $\Sigma =\{ 1 , \dots , N \}$ and transition probability matrix $P = [ p_{i j } ]$. The transition matrix operates from the right $x_{t+1} = x_t P$. A lumping  is defined by a partition of the state space $\widetilde{ \Sigma } = \{ L_1 , \dots , L_M \}$ where $L_k$ is a nonempty subset of $\Sigma$, $L _k \cap L_l = \O$ if $k \neq l$ and $\bigcup _k L_k = \Sigma$. By necessity $M \leq N$. A lumping can be defined by an $N \times M$ matrix $\Pi = [ \pi _{ik} ]$ where $\pi _{i k} =1$ if $i \in L_k$ and $\pi _{i k} =0$ if $i \notin L_k$. The lumping $\widetilde{ \Sigma }$ induces a quotient process $\tilde{x} _t = x_t \Pi$. If the process $\tilde{x}_t$ is a Markov process (which is not usually the case), then we say that the Markov chain is {\em strongly lumpable} with respect to $\Pi$ (or $\widetilde{ \Sigma }$). The following criterion is  necessary and sufficient for a Markov chain with transition matrix $P$ to be lumpable with respect to $\Pi$ or $\widetilde{ \Sigma }$ \cite{Kemeny99}: \\

\noindent {\rm THEOREM 1.} {\em If, for any $L_k , L_l \in \widetilde{ \Sigma }$, the total probability of a transition from any state $i \in L_k$ to $L_l$, i.e. $\sum _{j \in L_l} p_{ij}$, is independent of $i$, then $P$ is lumpable with respect to $\widetilde{\Sigma}$.} \\

\noindent The transition matrix for the reduced Markov chain is given by
\begin{equation}
	\tilde{p} _{k l } = \sum _{j \in L_l} p_{i j } \;\;\;\; i \in L_k ,
\label{p_reduced}
\end{equation}
where we note that $\tilde{P} = [ \tilde{p} _{kl} ]$ is well defined since the sum is independent of which representative $i \in L _k$ we chose. 

As mentioned in the introduction, Theorem 1 is not immediately useful for identifying lumpings of a Markov chain. In previous work by Barr and Thomas~\cite{Barr77}, a necessary lumpability criterion involving the left eigenvectors of the transition matrix was presented. It was first noted that the spectrum of $\tilde{P}$ must be a subset of the spectrum of $P$. Furthermore, it was noted that if $v P = \lambda v$ then $\tilde{v} \tilde{P} = \lambda \tilde{v}$, with $\tilde{v} = v \Pi$. It follows that if $\lambda$ is an eigenvalue of both $P$ and $\tilde{P}$, then $\tilde{v}$ is an eigenvector of $\tilde{P}$, but if $\lambda$ is not an eigenvalue of $\tilde{P}$ then $\tilde{v} = v \Pi =0$. Intuitively this observation can be understood as $\Pi$ eliminating a subset of the eigenvectors of $P$. It is an example of a general result concerning reduction of a linear operator $P$, which can be formulated as a commutation relation $\Pi \tilde{P}  =  P \Pi$ that must be fulfilled for a reduction $\Pi$ to be well defined~\cite{Lorch}. It is implied that $\mbox{ker} ( \Pi )$ is $P$-invariant (i.e. spanned by eigenvectors of $P$, if we for the moment ignore the problem of degenerated eigenvalues).

 Barr and Thomas'  result suggests a search for lumpings defined by $\Pi$ such that $v ^{\alpha} \Pi = 0$ for some subset of the left eigenvectors of $P$, $\{ v^{\alpha} \}_{\alpha \in J}$, $J \subseteq \Sigma$. Since $\Pi$ should be a matrix with zeros and ones, $v ^{\alpha} \Pi = 0$ essentially means searching for eigenvectors with subsets of elements that sums to zero. For lumpings only involving agglomeration of two states this is straightforward since the eigenvector(s) must have pairs of elements $v^{\alpha} _i = - v^{\alpha} _j$. However,  agglomeration of $k$ states means searching for partial sums evaluating to zero and involving $k$ elements. In addition (as we will see) there must be several eigenvectors simultaneously  satisfying $v^{\alpha} \Pi = 0$ (the criterion  is only necessary, not sufficient). As Barr and Thomas point out, there is no obvious algorithm to generate $\Pi$ based on their result.

In this paper we present a result that can easily be used for identifying possible partitions of the state space leading to a reduced Markov process. Our method is based on the observation that the dual of the probability space can be used to identify lumplings. The criterion $v^{\alpha} \Pi =0$ is viewed as an orthogonality condition between the left eigenvectors $\{ v^{\alpha} \}_{\alpha \in J}$ and the column space of $\Pi$. The orthogonal complement of a set of left eigenvectors  is spanned by complementary right eigenvectors (defined naturally in the dual vector space). These complementary eigenvectors span the column space of $\Pi$. Requiring that $\Pi$ consists of zeros and ones does, as we will see, correspond to a criterion of repeated elements within each complementary right eigenvector. Clearly, identifying repeated elements in the right eigenvectors is algorithmically straight forward. The precise result reads as follows. \\

\noindent {\rm THEOREM 2.} {\em Assume that $P$ is a diagonalizable transition matrix with full rank\footnote{The transition matrix $(1-\zeta ) P + \zeta 1$, $0 \leq \zeta < 1$, allows exactly the same lumping as $P$ according to Theorem 1. The rank condition is therefore not a real restriction.} describing a Markov process $x_{t+1} = x_t P$. Consider a set of linearly independent right eigenvectors of $P$, $P u^{\alpha} = \lambda ^{\alpha}  u^{\alpha}$. Let $I \subseteq \Sigma$ be the set of indices for the eigenvectors. Form state equivalence classes defined by states with identical elements in all eigenvectors $u^{\alpha}$, i.e. $i \sim j$ iff $u^{\alpha} _i =  u^{\alpha} _j$ $\forall \alpha \in I$. The equivalence classes  define a partitioning $\widetilde{ \Sigma }$ of the state space. This partitioning is a lumping of the Markov chain if the number of partition elements equals the number of eigenvectors, i.e. $| \widetilde{ \Sigma } | = | I |$. 

Conversely, if $\widetilde{\Sigma}$ is a lumping then there exist $| \widetilde{\Sigma} |$ linearly independent right eigenvectors that are invariant under permutations within the lumps.
} \\

\noindent {\em Proof.} Let $\{ v^{\beta} \}_{\beta = 1 , \dots , N}$  be the linearly independent ($1 \times N)$ left eigenvectors of $P$ and $\{ u^{\beta} \}_{\beta = 1 , \dots , N}$ be the ($N \times 1$) right eigenvectors. We normalize the left and/or right eigenvectors so that $v^{\beta}  u^{\gamma} = \delta ^{\beta \gamma}$ (where $\delta$ is the Kronecker delta). Let $\Pi$ define the partitioning $\widetilde{\Sigma} = \{ L _1 , \dots , L_{|I|} \}$. Since $\{ u^{\alpha} \} _{\alpha \in I}$ are linearly independent and $| \widetilde{ \Sigma } | = | I |$,  by the construction of the equivalence classes, the  column space of $\Pi$ is also spanned by $\{ u^{\alpha} \} _{\alpha \in I}$. From the orthogonality of left and right eigenvectors it follows that
\begin{equation}
	v^{\beta} \Pi = 0 , \;\;\;\; \forall \beta \notin I .
\label{v_orth_pi}
\end{equation}
We now decompose $P$ as
\begin{equation}
	p_{ij} = \sum _{\beta = 1}^N \lambda ^{\beta}  u ^{ \beta } _i v^{ \beta } _j , 
\label{P_decomposition}
\end{equation}
where $u^{\beta} _i$ and $v^{\beta} _i$ denote the $i$th entry in right respective left eigenvector $\beta$.  Let $i \in L_k$. From Eq.~\ref{P_decomposition} follows that
\begin{equation}
	\sum _{j \in L_l} p _{i j} = \sum _{\beta = 1}^N \lambda ^{ \beta }  u ^{ \beta } _i \sum _{j \in L_l} v^{ \beta } _j = \sum _{\beta =1} ^N \lambda ^{ \beta }  u ^{ \beta } _i \left( v ^{\beta} \Pi \right) _l = 
	\sum _{\beta \in I} \lambda ^{ \beta }  u ^{ \beta } _i  \left( v ^{\beta} \Pi \right) _l ,
\label{i_indep}
\end{equation}
where we use Eq.~\ref{v_orth_pi} in the last step. By the construction ($u_i ^{\beta} = u _j ^{\beta}$ if $i,j \in L_k$, $\beta \in I$)   it is clear that the right hand side of Eq.~\ref{i_indep} does not depend on $i \in L_k$. From Theorem 1 it follows that $\widetilde{ \Sigma}$ is a lumping of the Markov chain. 

We now show that if the Markov chain is lumpable with $ \widetilde{ \Sigma } $ defining the partition, then there exist $ | \widetilde{ \Sigma } | $ eigenvectors with constant elements over the lumps. Let  $\widetilde{ \Sigma } = \{ L _1 , \dots , L_{|I|} \}$ be a lumping of the Markov chain. Let $y$ be a vector defined as
\begin{equation}
	y _ i = \tilde{y} _ k \; \; \; i \in L_k ,
\label{u_ansatz}
\end{equation}
i.e. constant elements within the lumps defined by $\widetilde{ \Sigma }$, and $ \tilde{y}$ at this point an arbitrary vector on the reduced state space. Now consider
\begin{equation}
	\sum _j p_{i j} y_j = \sum _l \sum _{j \in L_l} p_{i j} y_j = \sum _l \tilde{p} _ {k l} \tilde{y} _l \;\;\; i \in L_k . 
\end{equation}
It follows that $y$ is a right eigenvector of $P$ if $\tilde{y}$ is an eigenvector of $\tilde{P}$. Since $\tilde{P}$ has exactly $| \widetilde{ \Sigma } |$ eigenvectors (clearly, $\tilde{P}$ is diagonalizable if $P$ is) it follows that $P$ has $| \widetilde{ \Sigma } |$ eigenvectors with the structure imposed in Eq.~\ref{u_ansatz}. $ \Box$

The two trivial lumpings: $\widetilde{ \Sigma } = \{ L _1 , \dots , L_N \}$ (all states in different lumps)  represented by all the right eigenvectors (since $P$ is full rank there can be no index pair equal in all eigenvectors); and $\widetilde{ \Sigma } = \{ L _1 \}$ (all states in one lump) represented by the always present $u = (1, 1 , \dots , 1 ) ^T$ eigenvector (implied by $\sum _j p_{ij} = 1$). As a consequence of $\widetilde{ \Sigma } =  \{ L _1 \}$ always being an accepted lumping, any eigenvector suggesting a unique two-lump partition is an accepted lumping. In practice, the right eigenvectors involved in forming the equivalence classes often, but not always, form a hierarchy of refined partitions (always including $\widetilde{ \Sigma } =  \{ L _1 \}$ as the most coarse partition), and the overall partition is set by the eigenvector suggesting the finest partition. It should be noted however that nothing prevents $\tilde{u}_k = \tilde{u}_l$ for $k\neq l$ in Eq.~\ref{u_ansatz}, which implies that not all  lumpings of a Markov chain must correspond to a single eigenvector with exactly that suggested lumping. A lumping can instead be composed of eigenvectors each by  themselves suggesting a coarser lumping. Consider for instance $(1, 1 , \dots , 1 ) ^T$, which is always a member of a subset of eigenvectors that correspond to a valid lumping. The four-lump partition ($\widetilde{ \Sigma} _5$) in Example 2 below also exemplifies this situation. In general, the total number of possible lumpings of a Markov chain may exceed the number of eigenvectors, i.e. the number of states. 

Theorem 2 is valid also if $P$ has redundant eigenvalues. However, eigenvectors corresponding to redundant eigenvalues are not unique. By "rotating"  the eigenvectors corresponding to a redundant eigenvalue it can be possible to find different suggested partitions. Especially, any $k$ elements can be set equal by choosing an appropriate linear combination of the $k$ eigenvectors corresponding to a $k$-fold degenerate eigenvalue. It can happen that forcing a subset of the elements to be equal simultaneously sets some other elements equal as well, creating a more coarse lumping than expected. The degenerate eigenvalues in Example 2 demonstrates this point. 

The identity matrix is an extreme example of redundant eigenvalues. In this case every partition of the state space is a lumping. The number of possible lumpings equals the Bell number $B_N$ and grows extremely fast with $N$. From this argument it is clear that, in some cases, generating all possible lumpings of a Markov chain is not computationally feasible. The lumpings suggested by the eigenvectors can be viewed as "generators" that can possibly be combined to create finer partitions of the state space. A more non-trivial example than the identity matrix is given in Example 2.  

\section{Examples}
Theorem 2 is illustrated with two examples. \\

\noindent {\bf Example 1.} Consider the transition matrix

\begin{equation*}
P = \left(
\begin{array}{ccc}
a +b+(c-1)/2 & 1-a-b & (1-c)/2   \\
-a+(c+1)/2 & a & (1-c)/2  \\ 
 1-b-c & b & c  
 \end{array}
\right) ,
\end{equation*}
with $0 \leq a,b,c \leq 1$.
$P$ has the eigenvalues $ \lambda ^1 = 1$, $\lambda ^2 = 2 a + b -1$ and $\lambda ^3 = (3c-1)/2$, and the eigenvectors  
\begin{eqnarray*}
	u^1 & = & (1 , 1 ,1  )\tp, \\  
	u^2 & = & ( 1+c-2a-2b , 2(a-c) , 2b+c-1 )\tp \text{  and}\\
	u^3 & = & ( -1 , -1 , 2 )\tp. 
\end{eqnarray*}
There are three possible lumpings of $P$:
\begin{equation*}
\begin{array}{lcll}
	\widetilde{ \Sigma}_1 & = &  \{ \{ 1, 2, 3 \} \} & \mbox{from $\{u^1\}$,} \\
	\widetilde{ \Sigma}_2 & = &  \{ \{ 1, 2 \}, \{ 3 \} \} & \mbox{from $\{u^1, u^3\}$ and} \\
	\widetilde{ \Sigma}_3 & = &  \{ \{ 1\}, \{ 2 \}, \{ 3 \} \} & \mbox{from $\{u^1, u^2, u^3\}$,}
\end{array}
\end{equation*}	
where $\widetilde{ \Sigma}_2$ is valid if $2a+b-1 \neq 0$ and $3c-1 \neq 0$. \\

\noindent {\bf Example 2.} Consider the transition matrix (from \cite{Cobb2003})

\begin{equation*}
P = \left(
\begin{array}{cccccccc}
0& a& a& a& a& 0& a& 0 \\
      a& 0& a& a& 0& a& 0& a\\
      a& a& 0& a& 0& a& a& 0\\
      a& a& a& 0& a& 0& 0& a\\
      b& 0& 0& b& 0& 0& b& b\\
      0& b& b& 0& 0& 0& b& b\\
      b& 0& b& 0& b& b& 0& 0\\
      0& b& 0& b& b& b& 0& 0
\end{array}
\right) ,
\end{equation*}
where $a=1/5$ and $b=1/4$. $P$ has two pairs of degenerate eigenvalues.
The associated (not normalized) right eigenvectors $u^i$ of $P$ have the following structure (columns sorted in order of decreasing $| \lambda |$):
\begin{equation*}
T =  (u^1 u^2 ... u^8) = \left(
\begin{array}{cccccccc}
	1& 0      & 0            & - 1       &  0           &  1         & -1       &  1    \\
         1& 0      &  0           & 1          & 0            &  -1        & -1       & 1    \\
         1& 0      & -1           &  0         &  -1         &  0          &  1       & 1    \\
         1& 0      &  1            &  0         &   1         &  0          &  1       & 1     \\
         1& - 1 &  - c_1        & c_1      &   c_2    &  c_2      & 0        & -5/4 \\
         1& - 1 &   c_1        & - c_1    &   - c_2   & -c_2      & 0       & -5/4 \\
         1& 1  &    c_1        & c_1       &     - c_2 &   c_2     & 0       & -5/4  \\
         1& 1  &  - c_1          & - c_1   &  c_2      & -c_2      & 0      & -5/4
\end{array}
\right) ,
\end{equation*}
for some numerical constants $c_1$ and $c_2$, ($c_1 \approx 0.579$ and $c_2 \approx 1.079$). Note that $u^3$ and $u^4$ share eigenvalue and therefore any linear combination $\alpha u^3 + \beta u^4$ is also an eigenvector. The same is true for $u^5$ and $u^6$.

The possible lumpings of $P$ and the corresponding subsets of right eigenvectors are given in Table \ref{Ex2Table}. Note that $\widetilde{\Sigma}_5$ is a lumping with an agglomeration not suggested by any of the individual eigenvectors involved, but rather by the set of eigenvectors as a whole. Further, partition $\widetilde{ \Sigma }_6$ through $\widetilde{ \Sigma }_9$ involve appropriate linear combinations of eigenvalues that correspond to the two degenerate eigenvalues of $P$, as indicated in Table~\ref{Ex2Table}.

\begin{table}
\caption{Possible lumpings of  $P$ and associated subsets of eigenvectors.}
\begin{center}
\begin{tabular}{|l|l|}
\hline
Lumping & Subset of eigenvectors \\
\hline
$\widetilde{ \Sigma } _1=\{ \{ 1, 2, 3, 4, 5, 6, 7, 8 \} \}$ & $\{u^1\}$\\
$\widetilde{ \Sigma } _2=\{ \{ 1, 2, 3, 4 \}, \{ 5, 6, 7, 8 \} \}$ & $\{u^1, u^8\}$\\
$\widetilde{ \Sigma } _3=\{ \{ 1, 2 \} , \{ 3, 4 \} , \{ 5, 6, 7, 8 \} \}$ & $\{u^1, u^7, u^8\}$\\
$\widetilde{ \Sigma } _4=\{ \{ 1, 2, 3, 4 \} , \{ 5, 6 \} , \{ 7, 8 \} \}$ & $\{u^1, u^2, u^8\}$\\
$\widetilde{ \Sigma }  _5  =  \{ \{ 1, 2 \} , \{ 3, 4 \} , \{ 5, 6 \} , \{ 7, 8 \} \}$ & $\{u^1, u^2, u^7, u^8\}$\\
$\widetilde{ \Sigma } _6  =  \{ \{ 1, 2 \} ,  \{ 3 \} ,  \{ 4 \} , \{ 5, 8 \},  \{ 6, 7 \} \}$ & $\{u^1,  u^3  , u^5  , u^7, u^8\}$  \\
$\widetilde{ \Sigma } _7  =  \{ \{ 1\} , \{ 2 \} , \{ 3, 4 \} ,  \{ 5, 7 \} , \{ 6, 8 \} \}$ & $\{u^1,  u^4 ,  u^6  , u^7, u^8\}$ \\
$\widetilde{ \Sigma } _8  =  \{ \{ 1, 3 \} , \{ 2, 4 \} ,  \{ 5, 6 \} , \{ 7 \} , \{ 8 \} \}$ & $\{u^1, u^2, u^3 +u^4  ,  u^5 - u^6  , u^8\}$ \\
$\widetilde{ \Sigma } _9  =  \{ \{ 1, 4 \} , \{ 2, 3 \} , \{ 5 \} , \{ 6 \} , \{ 7, 8 \} \}$ & $\{u^1, u^2,  u^3 - u^4  , u^5 + u^6  , u^8\}$\\
$\widetilde{ \Sigma } _{10}  =  \{ \{ 1 \},  \{ 2 \} , \dots , \{ 7 \} , \{ 8 \} \}$ & $\{u^1, u^2, ..., u^7, u^8\}$\\
\hline
\end{tabular}\end{center}
\label{Ex2Table}
\end{table}

\section{Discussion and conclusions}
We have presented a sufficient and necessary condition for Markov chain lumpability. The condition can be applied to identify a set of lumpings of a given Markov chain by identifying sets of identical elements in the dual eigenvectors. 
The results presented could have been expressed naturally using symmetries. The equivalence relation between identical elements in the dual eigenvectors implies invariance under a permutation symmetry. The lumpings are orbits of the permutation group, which acts as a normal subgroup extension to the transition matrix. In some situations this formulation can be advantageous, especially to clarify correspondence with general reduction methods in dynamical systems. In this paper we have chosen not to emphasize the symmetry formulation since it tends to obstruct the clarity of the presentation.

\bibliographystyle{siam}
\bibliography{markov_lumpability}

\begin{thebibliography}{10}

\bibitem{Barr77}
{\sc D.~R. Barr and M.~U. Thomas}, {\em An eigenvector condition for {M}arkov
  chain lumpability}, Operations Research, 25 (1977), pp.~1028--1031.

\bibitem{Buchholz}
{\sc P.~Buchholz}, {\em Hierarchical {M}arkovian models: symmetries and
  reduction}, Performance Evaluation, 22 (1995), pp.~93--110.

\bibitem{Cobb2003}
{\sc G.~W. Cobb and Y.~Chen}, {\em An application of {M}arkov chain {M}onte
  {C}arlo to community ecology}, The American Mathematical Monthly, 110 (2003),
  pp.~265--288.

\bibitem{Derisavi03}
{\sc S.~Derisavi, H.~Hermanns, and W.~H. Sanders}, {\em Optimal state-space
  lumping in {M}arkov chains}, Information Processing Letters, 87 (2003),
  pp.~309--315.

\bibitem{Gurvits}
{\sc L.~Gurvits and J.~Ledoux}, {\em Markov property for a function of a
  {M}arkov chain: A linear algebra approach}, Linear Algebra and its
  Applications, 404 (2005), pp.~85--117.

\bibitem{Jernigan03}
{\sc R.W. Jernigan and R.H. Baran}, {\em Testing lumpability in {M}arkov
  chains}, Statistics and Probability Letters, 64 (2003), pp.~17--23(7).

\bibitem{Kemeny99}
{\sc J.~G. Kemeny and J.~L. Snell}, {\em Finite {M}arkov Chains}, Springer, New
  York, NY, USA, 2nd~ed., 1976.

\bibitem{Lorch}
{\sc E.R. Lorch}, {\em Spectral theory}, Oxford University Press, New York,
  1962.

\bibitem{Obal}
{\sc W.D. Obal, M.G. McQuinn, and W.H. Sanders}, {\em Detecting and exploiting
  symmetry in discrete-state {M}arkov models}, in 12th Pacific Rim
  International Symposium on Dependable Computing (PRDC'06), Los Alamitos, CA,
  USA, 2006, IEEE Computer Society, pp.~26--38.

\bibitem{Ring}
{\sc A.~Ring}, {\em Symmetric lumping of {M}arkov chains}, in Jahrestagung der
  DMV, Halle, 2002, Deutsche Mathematiker-Vereinigung, p.~142.

\bibitem{Rogers}
{\sc L.C.G. Rogers and J.W. Pitman}, {\em Markov functions}, Annals of
  probability, 9 (1981), pp.~573--582.

\end{thebibliography}

\end{document}